\documentclass[12pt,a4paper,oneside,reqno,notitlepage]{amsart}
\usepackage{amssymb}
\textwidth
150mm

\theoremstyle{plain}

\numberwithin{equation}{section}

\begin{document}

\baselineskip 8mm
\parindent 9mm

\title[]
{Explosive solutions of parabolic stochastic partial differential equations with L$\acute{e}$vy noise}

\author{Kexue  Li,  Jigen Peng,  Junxiong Jia}

\address{Kexue Li\newline
School of Mathematics and Statistics, Xi'an Jiaotong University, Xi'an 710049, China; Department of Mathematics, University of Washington, Seattle, WA, 98195, USA}
\email{kexueli@gmail.com}

\address{Jigen Peng\newline
School of Mathematics and Statistics,
Xi'an Jiaotong University, Xi'an 710049, China}
\email{jgpeng@mail.xjtu.edu.cn}

\address{Junxiong Jia\newline
School of Mathematics and Statistics,
Xi'an Jiaotong University, Xi'an 710049, China}
\email{jiajunxiong@163.com}

\thanks{{\it 2010 Mathematics Subjects Classification}: 60H15; 60J75}

\begin{abstract}
In this paper, we study the explosive solutions to a class of parbolic stochastic semilinear differential equations driven by a L$\acute{\mbox{e}}$vy type noise. The sufficient conditions are presented to guarantee the existence of a unique positive solution of the stochastic partial differential equation under investigation. Moreover, we show that positive solutions will blow up in finite time in mean $L^{p}$-norm sense, provided that the initial data, the nonlinear term and the multiplicative noise satisfies some conditions. Several examples are presented to illustrated the theory. Finally, we establish a global existence theorem based on a Lyapunov functional and prove that a stochastic Allen-Cahn  equation driven by L$\acute{\mbox{e}}$vy noise has a global solution.
\end{abstract}
\maketitle

\section{\textbf{Introduction}}
Fujita \cite{Fujita}  considered the initial-boundary problem for a semilinear parabolic equation
\begin{equation}\label{Fu}
\left\{\begin{aligned}
&\frac{\partial u}{\partial t}=\Delta u+u^{1+\alpha}, \ t>0, \ x\in \mathbb{R}^{d}, \\
&u(x,0)=a(x),\ x\in \mathbb{R}^{d},\\
\end{aligned}\right.
\end{equation}
Fujita  showed that there does not exist a global solution for any nontrivial nonnegative initial data when $0<d\alpha<2$,  and there exists a global solution for sufficiently small initial data when $d\alpha>2$. Hayakawa \cite{Hayakawa} proved that (\ref{Fu}) has no global solution for any nontrivial nonnegative initial data in the critical case $d\alpha=2$ if the dimension $d$ equals 1 or 2. \\
Fujita \cite{Fu} studied the initial-boundary problem for a semilinear parabolic equation in domain $D\subset \mathbb{R}^{d}$:
\begin{equation}\label{Fujita}
\left\{\begin{aligned}
&\frac{\partial u}{\partial t}=\Delta u+f(u), \ t>0, \ x\in D, \\
&u(x,0)=a(x),\ x\in D,\\
&u(x,t)=0, \ x\in \partial D,
\end{aligned}\right.
\end{equation}
Denote $\lambda_{0}$ as the smallest eigenvalue of $-\Delta$ and the corresponding eigenfunction $\phi_{0}>0$, satisfies $\int_{D}\phi_{0}(x)dx=1$ and
\begin{equation*}
\left\{\begin{aligned}
&-\Delta \phi_{0}=\lambda_{0}\phi_{0}  \ \ \  \mbox{in} \ D,\\
&\phi_{0}|_{\partial D}=0.
\end{aligned}\right.
\end{equation*}
Assume that $f$ satisfies the following \\
(\mbox{f.1}) $f$ is locally Lipschitz continuous. \\
(\mbox{f.2}) $f(0)\geq 0$  and $f(r)>0$ for $r>0$. \\
(\mbox{f.3}) $1/f$ is integrable at $t=+\infty$. \\
(\mbox{f.4}) $f$ is convex in $[0,\infty)$. \\
(\mbox{f.5}) $f(r)-\lambda_{0}r>0 \ \mbox{for} \ \ r>\int_{D}a_{0}\phi_{0}dx$, where $a_{0}=\exp(-k|x|^{2})$, $k>0$, $x\in \mathbb{R}^{d}$.\\
Fujita \cite{Fu} showed that if $D$ is bounded, $a(x)\geq 0$ in $D$ and $f$ satisfies (\mbox{f.1})- (\mbox{f.5}), then the solution of  (\ref{Fujita})  blows up in finite time. As a special case that $f(r)=r^{1+\alpha}(\alpha> 0)$, the solution of (\ref{Fujita}) blows up in finite time if
\begin{equation*}
\int_{D}a(x)\phi_{0}(x)dx\geq \lambda^{1/\alpha}_{0}.
\end{equation*}
We refer to \cite{GV} about the many developments on solutions of nonlinear parabolic equations may blow up in finite time.

Recent years, stochastic partial differential equations has attracted the attention of many researchers. It is of interest to study the non-existence of global solutions to parabolic stochastic partial differential equations perturbed by random noise as follows:
\begin{equation}\label{multiplicative}
\left\{\begin{aligned}
&du=\Delta u+f(u)+\sigma(u)dW_{t}, \ t>0, \ x\in D, \\
&u(x,0)=g(x),\ x\in D,\\
&u(x,t)=0, \ x\in \partial D.
\end{aligned}\right.
\end{equation}
When $f(u)\equiv0$, $\sigma(u)=u^{\gamma}\ (\gamma\geq 1)$,
Mueller \cite{Mueller}  considered the equation
\begin{equation}\label{gamma}
\left\{\begin{aligned}
&\frac{\partial u}{\partial t}=\Delta u+u^{\gamma}\dot{W}, \ \gamma\geq 1, \ t>0, \ 0\leq x\leq J,\\
&u(t,0)=u(t,J)=0,
\end{aligned}\right.
\end{equation}
where $\dot{W}=\dot{W}(t,x)$ is 2-parameter white noise and $u(x,0)$ is nonnegative and continuous. The conclusion is that for $1\leq \gamma<\frac{3}{2}$, $u$ exists for all time. Mueller \cite{M} showed that when $u(x,0)$ is a continuous nonnegative function on $[0,J]$, vanishing at the endpoints, but not identically zero, then there is a positive probability that the solution $u$ of (\ref{gamma}) blows up in finite time if $\gamma>3/2$. When $\sigma(u)\equiv 1$ and the Laplacian operator $\Delta$ is replaced by the infinitesimal generator of a $C_{0}$ semigroup,  Prato and Zabczyk \cite{PZ} considered the stochastic semilinear equation
\begin{equation}\label{prato}
\left\{\begin{aligned}
&du=(Au+F(u))dt+dW\\
&u(0)=\xi,
\end{aligned}\right.
\end{equation}
where $A$ is the generator of semmigroup $S(t)=e^{At}$ on a Banach space $E$, and $F$ is a mapping from $E$ into  $E$. $W$ is a Wiener process defined on a probability space $(\Omega, \mathcal{F}_{t}, P)$. $\xi$ is an $\mathcal{F}_{0}$-measurable $E$-valued random variable.  They assumed that $F$ satisfies the Lipschitz condition on bounded sets of  $E$.  This property of $F$ together with some other conditions ensure that (\ref{prato}) has a unique non-exploding solution. When $\sigma(u)=\sigma$ (positive constant), $W$ is a 2-dimensional Brownian sheet, $f$ is a nonnegative, convex function such that $\int_{0}^{\infty}1/f<\infty$, Bonder and Groisman \cite{BG} proved that the solution to (\ref{multiplicative}) blows up in finite time with probability one for every nonnegative initial datum $u(x,0)\geq 0$. Dozzi and L$\acute{\mbox{o}}$pez-Mimbela \cite{DL} considered the equation (\ref{multiplicative}) with $\sigma(u)=\kappa u$,  where  $\beta, \kappa$ are given positive numbers, $W_{t}$ is a standard one-dimensional Brownian motion.  They proved that the solution of (\ref{multiplicative}) blows up in finite time with positive probability  if $f(u)\geq Cu^{1+\beta}$ ($u>0$, $C>0$ ).

Chow \cite{Chow,Ch} considered the initial-boundary value problem for the parabolic It$\hat{\mbox{o}}$ equation
\begin{equation}\label{Ito}
\left\{\begin{aligned}
&\frac{\partial u}{\partial t}=Au+f(u,x,t)+\sigma(u,\nabla u,x,t)\partial_{t}W(x,t)\\
&u(x,0)=g(x),\ x\in D,\\
&u(x,t)=0, \ t\in (0,T),x\in \partial D,
\end{aligned}\right.
\end{equation}
where $D\subset \mathbb{R}^{d}$, $A=\sum_{i,j=1}^{d}\frac{\partial}{\partial x_{i}}[a_{ij}(x)\frac{\partial}{\partial x_{j}}]$ is a symmetric, uniformly elliptic operator with smooth coefficients, $f$ and $\sigma$ are given functions.  For $x\in \mathbb{R}^{d}, \ t\geq 0$, $W(x,t)$ is a continuous Wiener random field defined in a complete probability space  $(\Omega, \mathcal{F}, P)$ with a filtration $\mathcal{F}_{t}$. $W(x,t)$ has mean $\mathbb{E}W(x,t)=0$ and covariance function $q(x,y)$ defined by
\begin{equation*}
\mathbb{E}\{W(x,t)W(y,t)\}=(t \wedge s)q(x,y), \ x,y\in \mathbb{R}^{d},
\end{equation*}
where $t \wedge s=\min\{t,s\}$ for $0\leq t,s\leq T$.  Let $H=L^{2}(D)$, $H^{1}=H^{1}(D)$ be the $L^{2}$-Sobolev space of first order and $H^{1}_{0}$ the closure in $H^{1}$ of the space of $C^{1}$-functions with compact support in $D$. Under the usual conditions, such as coercivity conditions, Lipschitz continuity and boundedness conditions, Eq.(\ref{Ito}) has a a unique global strong solution $u\in C([0,T];H)\cap L^{2}([0,T]; H^{1})$ (see \cite[Theorem 3-7.2]{C}). To consider positive (nonnegative) solutions, the author assume that the following conditions hold: \\
$\mbox{(P1)}$ There exists a constant $\delta\geq 0$ such that
\begin{equation*}
\frac{1}{2}q(x,x)\sigma^{2}(r,\xi,x,t)-\sum_{i,j=1}^{d}a_{ij}(x)\xi_{i}\xi_{j}\leq \delta r^{2},
\end{equation*}
for all $r\in \mathbb{R}$, $x\in \overline{D}$, $\xi\in \mathbb{R}^{d}$ and $t\in [0,T]$.\\
$\mbox{(P2)}$ The function $f(r,x,t)$ is continuous on $\mathbb{R}\times \overline{D}\times [0,T]$ such that $f(r,x,t)\geq 0$ for $r\leq 0$ and $x\in \overline{D}$, $t\in [0,T]$.\\
$\mbox{(P3)}$ The initial datum $g(x)$ on $\overline{D}$ is positive and continuous. \\
Chow \cite{Chow} proved that the solution of Eq.(\ref{Ito}) is positive. Under some suitable conditions, Chow \cite{Chow, Ch} showed that the positive solutions of a class of stochastic reaction-diffusion equations will blow up in the $L^{p}$-norm sense, $p\geq 1$. Chow and Liu \cite{CL} considered the problem of explosive solutions in mean $L^{p}$-norm sense of semilinear stochastic functional parabolic  differential equations of retarded type.

Lv and Duan \cite{Lv} considered the Eq.(\ref{Ito}) with $A=\triangle$, the Laplacian operator, the nonlinear term $f$ is assumed to be satisfied by an inequality, which is weaker than the condition $\mbox{(P2)}$, the noise intensity $\sigma$ allows to be higher nonlinear than the square nonlinear (see \cite[formula (3.2)]{Lv}). They proved that the noise could induce finite time blow up of solutions.

Recent years, stochastic partial differential equations driven by L$\acute{\mbox{e}}$vy noise have attracted many attentions (see, for example,  \cite{ZJ, W, Y, L, PZ, D, LZ, RZ,  SH}). But there are few results about the existence of explosive solutions to stochastic partial differential equations with L$\acute{\mbox{e}}$vy noise  in the literature.  Bao and yuan \cite{Bao} considered the stochastic reaction-diffusion equations with jumps
\begin{equation}\label{jump}
\left\{\begin{aligned}
&\frac{\partial u}{\partial t}=Au+b(t,u,x)+\sigma(t,u,\nabla u,x)\partial_{t}W(x,t)\\
&\quad\quad\quad+\int_{\mathbb{Y}}\Upsilon(t,u,x,y)\partial_{t}\widetilde{N}(t,dy), \ t>0,\\
&u(x,0)=\phi(x),\ x\in \mathcal{O}, \ u(t,x)|_{\partial \mathcal{O}}=0, \ t>0,
\end{aligned}\right.
\end{equation}
where $\mathcal{O}\subset \mathbb{R}^{n}$ is a bounded domain with $C^{\infty}$ boundary $\partial \mathcal{O}$, $A=\sum_{i,j=1}^{n}\frac{\partial}{\partial x_{i}}(a_{i,j}(x)\frac{\partial}{\partial x_{j}})$ be a symmetric, uniformly elliptic operator with smooth coefficients, $W(x,t)$ is a Wiener random defined on the completed probability space $(\Omega,\mathcal{F},\{\mathcal{F}\}_{t\geq 0}, \mathbb{P})$, $\widetilde{N}(dt,du)$ is the compensated Poisson measure. Under some conditions, they showed that the solution of (\ref{jump}) blows up in finite time. It should be pointed out that the nonlinear term $b:[0,\infty)\times \mathbb{R}\times \overline{O}\mapsto \mathbb{R}$ is assumed to be locally Lip-continuous w.r.t the second variable such that $b(t,r,x)\geq 0$ for any $r\leq 0$, however, there are many functions don't satisfy this condition, for example, $b(r)=r(1-r^{2})$. And $\mathcal{O}\subset \mathbb{R}^{n}$ is assumed to be bounded, the proof of Theorem 2.1 in \cite{Bao} depends on the boundedness of volume of $\mathcal{O}$. The results of \cite{Bao} can't be generalized to the case for unbounded domain, such as $\mathcal{O}=\mathbb{R}^{n}$.

In this paper, we study the problem of explosive solutions to a class of semilinear stochastic parabolic differential equations driven by L$\acute{\mbox{e}}$vy noise. The paper is organized as follows. In Section 2, we recall some basic results for semilinear stochastic parabolic equations with L$\acute{\mbox{e}}$vy noise. In Section 3, under some assumptions, we prove that the existence of positive solutions of a semilinear stochastic reaction-diffusion equation. In Section 4, under some suitable conditions on the drift or diffusion term, we prove that the solutions of stochastic parabolic differential equations will blow up in a finite time in mean $L^{p}$-norm sense, $p\geq 1$. Some examples are presented to illustrate the theory. In Section 5, we establish a global existence theorem based on a Lyapunov functional. We show that the existence of global solution to stochastic Allen-Cahn equation driven by L$\acute{\mbox{e}}$vy noise.

\section{\textbf{Preliminaries}}

Let $D$ be a domain in $\mathbb{R}^{d}$, which has a smooth boundary if it is bounded. Denote $L^{2}(D)$ by $H$, the usual $L^{2}$ real Hilbert space with the inner product $(\cdot,\cdot)$ and norm $\|\cdot\|$, respectively. Let $H^{1}=H^{1}(D)$ be the $L^{2}$-Sobolev space of first order. Denote $H_{0}^{1}$ the closure in $H^{1}$ of the space of $C^{1}$-functions with compact support in $D$. Denote by $D([0,T],H)$ the space of all c$\grave{a}$dl$\grave{a}$g paths from $[0,T]$ into $H$.
Let $W(x,t)$ be a continuous Wiener random field defined on a complete probability space $(\Omega, \mathcal{F}, \mathbb{P})$ with a filtration $\mathcal{F}_{t}$. $W(x,t)$ has mean zero and covariance function $q(x,y)$ such that
\begin{equation*}
\mathbb{E}W(x,t)=0, \ \ \mathbb{E}\{W(x,t)W(y,t)\}=(t\wedge s)q(x,y),\ s,t\in [0,T], x,y\in \mathbb{R}^{d}.
\end{equation*}
The associated covariance operator $Q$ in $H$ with kernel $q(x,y)$ is defined by
\begin{equation*}
(Q\phi)(x)=\int_{D}q(x,y)\phi(y)dy, \ \ x\in D,\ \phi\in H.
\end{equation*}
In this paper, we assume that the covariance function $q(x,y)$ is bounded, continuous and there is $q_{0}>0$ such that
\begin{equation*}
\sup_{x,y\in D}|q(x,y)|\leq q_{0}\ \ \mbox{and} \ \mbox{Tr}\ Q=\int_{D}q(x,x)dx<\infty.
\end{equation*}
Let $(Z, \mathcal{B}(Z)$ be a  measurable space. Denote by $N(dt,dz)$ the Poisson random measure with intensity measure $dt\nu(dz)$ on $\mathbb{R}_{+}\times Z$, here $\mathbb{R}_{+}=[0,\infty)$, $dt$ is the Lebesgue measure on $\mathbb{R}_{+}$, $\nu(dz)$ is a $\sigma$-finite measure on $(Z, \mathcal{B}(Z))$. Denote by $\widetilde{N}(dt,dz)=N(dt,dz)-dt\nu(dz)$ the compensated Poisson measure. Assume that $W$ and $N$ are independent.

Consider the initial-boundary problem of a semilinear stochastic reaction-diffusion equation in domain $D\subset \mathbb{R}^{d}$:
\begin{equation}\label{rd}
\left\{\begin{aligned}
&\frac{\partial u}{\partial t}=Au+f(u,x,t)+\sigma(u,\nabla u,x,t)\partial_{t}W(x,t)\\
&\quad\quad+\int_{Z}\varphi(u,x,z,t)\partial_{t}\widetilde{N}(t,dz),\\
&u(x,0)=g(x),\ x\in D,\\
&u(x,t)=0, \ t\in (0,T),x\in \partial D,
\end{aligned}\right.
\end{equation}
where $A=\sum_{i,j=1}^{d}\frac{\partial}{\partial x_{i}}[a_{ij}(x)\frac{\partial}{\partial x_{j}}]$ is a symmetric, uniformly elliptic operator with smooth coefficients, that is, there exists a constant $c>0$ such that $b(x,\xi):=\sum_{i,j=1}^{d}a_{ij}(x)\xi_{i}\xi_{j}\geq c|\xi|^{2}$ for all $x\in \overline{D}$ and $\xi=(\xi_{1},\ldots, \xi_{d})\in \mathbb{R}^{d}$.

Let $u_{t}=u(\cdot,t)$, $F_{t}(u)=f(u,\cdot,t)$, $\Sigma_{t}(u)=\sigma(u,\nabla u, \cdot,t)$,  \ $\Gamma_{t}(u,z)=\varphi(u,\cdot,z,t)$ and $W_{t}=W(\cdot,t)$, then we can rewrite the equation (\ref{rd}) as
\begin{equation}\label{parabolic}
\left\{\begin{aligned}
&du_{t}=[Au_{t}+F_{t}(u_{t})]dt+\Sigma_{t}(u_{t})dW_{t}+\int_{Z}\Gamma_{t}(u_{t},z)\widetilde{N}(dt,dz),\\
&u_{0}=g,
\end{aligned}\right.
\end{equation}
where $A$ is regarded as a linear operator from $H^{1}$ into $H^{-1}$ with domain $D(A)=H_{0}^{1}\cap H^{2}$, $F_{t} :H\rightarrow H$ is continuous.
If $A$ satisfies the coercivity condition, $f$ and $\sigma$ satisfy the Lipschitz continuity and boundedness conditions, the equation (\ref{parabolic}) has a unique global strong solution $u\in L^{2}([0,T]; H^{-1})\cap D([0,T]; H)$( see Theorem 3.2, \cite{RZ}).

To consider the positive solutions, we assume that (\ref{rd}) has a unique (strong) solution. In addition, we assume that

$\mbox{(A1)}$
\begin{equation*}
f(u,x,t)\geq a_{1}u^{\beta}+a_{2}u,
\end{equation*}
where $a_{1},a_{2}\in \mathbb{R}$, $\beta>1$, $(-1)^{\beta}\in \mathbb{R}$ and
\begin{equation}
a_{1}\left\{\begin{aligned}
&>0, \ \mbox{if}\ (-1)^{\beta}=1,\\
&<0, \ \mbox{if}\ (-1)^{\beta}=-1.
\end{aligned}\right.
\end{equation}

$\mbox{(A2)}$ There exist constants $b_{1},b_{2}\geq 0$ such that
\begin{equation*}
\frac{1}{2}q(x,x)\sigma^{2}(u,\xi,x,t)-\sum_{i,j=1}^{d}a_{ij}(x)\xi_{i}\xi_{j}\leq b_{1}|u|^{m}+b_{2}u^{2},
\end{equation*}
for all $u\in \mathbb{R}$, $x\in \overline{D}$, $\xi\in \mathbb{R}^{d}$ and $t\in [0,T]$, where $2< m<\beta+1$.

$\mbox{(A3)}$ There exist a a constant $\mu\in [2,\beta+1)$ and mappings $\psi: \overline{D}\rightarrow \mathbb{R}_{+}$ with $\int_{Z}\psi(z)\nu(dz)<+\infty$, such that
\begin{equation}
\varphi^{2}(u,x,z,s)\leq \psi(z)|u(x,s)|^{\mu}.
\end{equation}

$\mbox{(A4)}$ the initial datum $g(x)$ on  $\overline{D}$ is positive and continuous.

As in \cite{Chow}, let $\eta(r)=r^{-}$ denote the negative part of $r$ for $r\in \mathbb{R}$, or $\eta(r)=0$, if $r\geq 0$ and
$\eta(r)=-r$ if $r<0$. Set $k(r)=\eta^{2}(r)$ so that $k(r)=0$ for $r\geq 0$ and $k(r)=r^{2}$ for $r<0$. For $\varepsilon>0$, let $k_{\varepsilon}(r)$ be a $C^{2}$-regularization of $k(r)$ defined by
\begin{equation}\label{kr}
k_{\varepsilon}(r)=\left\{\begin{aligned}
&r^{2}-\frac{\varepsilon^{2}}{6}, &r<-\varepsilon,\\
&-\frac{r^{3}}{\varepsilon}(\frac{r}{2\varepsilon}+\frac{4}{3}), \ &-\varepsilon\leq r<0,\\
&0, &r\geq 0.
\end{aligned}\right.
\end{equation}
It is easy to see that $k_{\varepsilon}(r)$ has the following properties. \\
\textbf{Lemma 2.1.} (see \cite{Chow}). The first two derivatives $k'_{\varepsilon}$,  $k''_{\varepsilon}$ of $k_{\varepsilon}$ are continuous and satisfy the conditions: $k'_{\varepsilon}(r)=0$ for $r\geq 0$; and $k''_{\varepsilon}(r)\geq 0$ for any $r\in \mathbb{R}$. Moreover, as $\varepsilon\rightarrow 0$, we have
\begin{align}\label{convergence}
k_{\varepsilon}(r)\rightarrow k(r),\ k'_{\varepsilon}(r)\rightarrow -2\eta(r)\ \  \mbox{and} \ k''_{\varepsilon}(r)\rightarrow 2\theta(r),
\end{align}
where $\theta(r)=0$ for $r\geq 0$, $\theta(r)=1$ for $r<0$, and the convergence is uniform for $r\in \mathbb{R}$.
\section{\textbf{Positive solutions}}
In this section, we will consider the existence of positive solution of Eq. (\ref{rd}).
\textbf{Theorem 3.1} Suppose that the conditionss $\mbox{(A1)-(A4)}$ hold. Then the solution of initial-boundary value problem (\ref{rd}) with nonnegative and continuous data remains positive so that $u(x,t)\geq 0$, a.s. for almost every $x\in D$ and for all $t\in [0,T]$. \\
\textbf{Proof.} Let $u_{t}=u(\cdot,t)$ and
\begin{equation*}
\Phi_{\varepsilon}(u_{t})=(1,k_{\varepsilon}(u_{t}))=\int_{D}k_{\varepsilon}(u(x,t))dx.
\end{equation*}
From It$\hat{\mbox{o}}$'s formula, it follows that
\begin{align*}
\Phi_{\varepsilon}(u_{t})&=\Phi_{\varepsilon}(g)+\int_{0}^{t}\int_{D}k'_{\varepsilon}(u(x,s))A u(x,s)dxds\\
&\quad+\int_{0}^{t}\int_{D}k'_{\varepsilon}(u(x,s))f(u(x,s),x,s)dxds\\
&\quad+\int_{0}^{t}\int_{D}k'_{\varepsilon}(u(x,s))\sigma(u(x,s),\nabla u(x,s),x,s)dW(x,s)dx\\
&\quad+\frac{1}{2}\int_{0}^{t}\int_{D}k''_{\varepsilon}(u(x,s))q(x,x)\sigma^{2}(u(x,s),\nabla u(x,s),x,s)dxds\\
&\quad+\int_{0}^{t}\int_{Z}\int_{D}\big(k_{\varepsilon}(u(x,s)+\varphi(u,x,z,s))-k_{\varepsilon}(u(x,s))\big)dx\widetilde{N}(dsdz)\\
&\quad+\int_{0}^{t}\int_{Z}\int_{D}\big(k_{\varepsilon}(u(x,s)+\varphi(u,x,z,s))-k_{\varepsilon}(u(x,s))-\varphi(u,x,z,s)k'_{\varepsilon}(u(x,s)))dx\nu(dz)ds\\
&\quad=\Phi_{\varepsilon}(g)+\int_{0}^{t}\int_{D}k''_{\varepsilon}(u(x,s))\big(\frac{1}{2}q(x,x)\sigma^{2}(u(x,s),\nabla u(x,s),x,s)-b(x,\nabla u(x,s))\big)dxds\\
&\quad+\int_{0}^{t}\int_{D}k'_{\varepsilon}(u(x,s))f(u(x,s),x,s)dxds\\
&\quad+\int_{0}^{t}\int_{D}k'_{\varepsilon}(u(x,s))\sigma(u(x,s),\nabla u(x,s),x,s)dW(x,s)dx\\
&\quad+\int_{0}^{t}\int_{Z}\int_{D}\big(k_{\varepsilon}(u(x,s)+\varphi(u,x,z,s))-k_{\varepsilon}(u(x,s))\big)dx\widetilde{N}(dsdz)\\
&\quad+\int_{0}^{t}\int_{Z}\int_{D}\big(k_{\varepsilon}(u(x,s)+\varphi(u,x,z,s))-k_{\varepsilon}(u(x,s))-\varphi(u,x,z,s)k'_{\varepsilon}(u(x,s)))dx\nu(dz)ds.
\end{align*}
By taking expectations of both sides of the above equality, we have
\begin{align*}
\mathbb{E}\Phi_{\varepsilon}(u_{t})&=\Phi_{\varepsilon}(g)+\mathbb{E}\int_{0}^{t}\int_{D}k''_{\varepsilon}(u(x,s))\big(\frac{1}{2}q(x,x)\sigma^{2}(u(x,s),\nabla u(x,s),x,s)-b(x,\nabla u(x,s))\big)dxds\\
&\quad+\mathbb{E}\int_{0}^{t}\int_{D}k'_{\varepsilon}(u(x,s))f(u(x,s),x,s)dxds\\
&\quad+\mathbb{E}\int_{0}^{t}\int_{Z}\int_{D}\big(k_{\varepsilon}(u(x,s)+\varphi(u,x,z,s))-k_{\varepsilon}(u(x,s))-\varphi(u,x,z,s)k'_{\varepsilon}(u(x,s)))dx\nu(dz)ds.
\end{align*}
From $\mbox{(A1)}$ and Lemma 2.1, it follows that
\begin{align}\label{inequality}
\mathbb{E}\Phi_{\varepsilon}(u_{t})&\leq \Phi_{\varepsilon}(g)+\mathbb{E}\int_{0}^{t}\int_{D}k''_{\varepsilon}(u(x,s))\big(b_{1}|u(x,s)|^{m}+b_{2}|u(x,s)|^{2}\big)dxds\nonumber\\
&\quad+\mathbb{E}\int_{0}^{t}\int_{D}k'_{\varepsilon}(u(x,s))(a_{1}u^{\beta}(x,s)+a_{2}u(x,s))dxds\nonumber\\
&\quad+\mathbb{E}\int_{0}^{t}\int_{Z}\int_{D}\big(k_{\varepsilon}(u(x,s)+\varphi(u,x,z,s))-k_{\varepsilon}(u(x,s))-\varphi(u,x,z,s)k'_{\varepsilon}(u(x,s)))dx\nu(dz)ds.
\end{align}
By Taylor's theorem, in view of the integral form of the remainder, we have
\begin{align}\label{Taylor}
k_{\varepsilon}(u(x,s)+\varphi(u,x,z,s))-k_{\varepsilon}(u(x,s))-(k'_{\varepsilon}(u(x,s)),\varphi(u,x,z,s))\nonumber\\
=\int_{0}^{1}(1-\tau)k''_{\varepsilon}(\varphi(u,x,z,s)\tau+u(x,s))\varphi^{2}(u,x,z,s)d\tau
\end{align}
Substitute (\ref{Taylor}) into (\ref{inequality}), we get
\begin{align}\label{integral}
\mathbb{E}\Phi_{\varepsilon}(u_{t})&\leq \Phi_{\varepsilon}(g)+\mathbb{E}\int_{0}^{t}\int_{D}k''_{\varepsilon}(u(x,s))\big(b_{1}|u(x,s)|^{m}+b_{2}|u(x,s)|^{2}\big)dxds\nonumber\\
&\quad+\mathbb{E}\int_{0}^{t}\int_{D}k'_{\varepsilon}(u(x,s))(a_{1}u^{\beta}(x,s)+a_{2}u(x,s))dxds\nonumber\\
&\quad+\mathbb{E}\int_{0}^{t}\int_{Z}\int_{D}\int_{0}^{1}(1-\tau)k''_{\varepsilon}(\varphi(u,x,z,s)\tau+u(x,s))\varphi^{2}(u,x,z,s)d\tau dx\nu(dz)ds.
\end{align}
Since $\lim_{\varepsilon\rightarrow 0}\mathbb{E}\Phi_{\varepsilon}(u_{t})=\mathbb{E}\|\eta(u_{t})\|^{2}$, taking the limits on both sides of (\ref{integral}) as $\varepsilon\rightarrow 0$, by (\ref{convergence}) we obtain
\begin{align}\label{limit}
\mathbb{E}\|\eta(u_{t})\|^{2}&\leq \int_{D}|\eta(g(x))|^{2}dx+2\mathbb{E}\int_{0}^{t}\int_{D}\theta(u(x,s))\big(b_{1}|u(x,s)|^{m}+b_{2}|u(x,s)|^{2}\big)dxds\nonumber\\
&\quad-2\mathbb{E}\int_{0}^{t}\int_{D}\eta(u(x,s))(a_{1}u^{\beta}(x,s)+a_{2}u(x,s))dxds\nonumber\\
&\quad+2\mathbb{E}\int_{0}^{t}\int_{Z}\int_{D}\int_{0}^{1}(1-\tau)\theta(\varphi(u,x,z,s)\tau+u(x,s))\varphi^{2}(u,x,z,s)d\tau dx\nu(dz)ds.
\end{align}
By the definition of $\eta$, it follows that $\eta(g)=0$. This together with $\mbox{(A4)}$, Lemma 2.1 and $(-1)^{\beta}a_{1}=|a_{1}|$ yield
\begin{align}\label{u}
\mathbb{E}\|\eta(u_{t})\|^{2}
&\leq 2\mathbb{E}\int_{0}^{t}\int_{D}\big(b_{1}|u(x,s)|^{m}+b_{2}|u(x,s)|^{2}\big)dxds\nonumber\\
&\quad-2\mathbb{E}\int_{0}^{t}\int_{D}\eta(u(x,s))(a_{1}u^{\beta}(x,s)+a_{2}u(x,s))dxds\nonumber\\
&\quad+\int_{Z}\psi(z)\nu(dz)\ \mathbb{E}\int_{0}^{t}\int_{D}|u(x,s)|^{\mu}dxds\nonumber\\
&= 2\mathbb{E}\int_{0}^{t}\int_{D}[b_{1}{(u^{-})}^{m}(x,s)+b_{2}{(u^{-})}^{2}(x,s)]dxds\nonumber\\
&\quad-2\mathbb{E}\int_{0}^{t}\int_{D}[|a_{1}|{(u^{-})}^{\beta+1}(x,s)-a_{2}{(u^{-})}^{2}(x,s)]dxds\nonumber\\
&\quad+\int_{Z}\psi(z)\nu(dz)\ \mathbb{E}\int_{0}^{t}\int_{D}{(u^{-})}^{\mu}(x,s)dxds.
\end{align}
It is known that the following $L^{p}$ interpolation inequality and Young inequality hold (see \cite{GT})
\begin{align}\label{alpha}
\|u\|_{L^{r}}\leq \|u\|^{\alpha}_{L^{p}}\|u\|^{1-\alpha}_{L^{q}},
\end{align}
\begin{align}\label{rpq}
ab\leq \varepsilon a^{\delta}+\varepsilon ^{-\frac{\omega}{\delta}}b^{\omega},
\end{align}
where $\alpha\in (0,1)$, $\varepsilon>0$, $\delta>0$, $\omega>0$, $a>0$, $b>0$,
\begin{align*}
&\frac{1}{r}=\frac{\alpha}{p}+\frac{1-\alpha}{q}, \ 1\leq p \leq r\leq q \leq \infty,\\
&\frac{1}{\delta}+\frac{1}{\omega}=1.
\end{align*}
Since $2< m<\beta+1$,  it follows that from (\ref{alpha}) and (\ref{rpq})
\begin{align}\label{beta}
2b_{1}\int_{D}(u^{-})^{m}(x,t)dx&=2b_{1}\|u^{-}\|^{m}_{L^{m}}\nonumber\\
&\leq C\|u^{-}\|^{m\alpha}_{L^{2}}\|u^{-}\|^{m(1-\alpha)}_{L^{\beta+1}}\nonumber\\
&\leq  \varepsilon \|u^{-}\|_{L^{\beta+1}}^{m(1-\alpha)\frac{2}{2-m\alpha}}+C(\varepsilon, m, \beta)\|u^{-}\|^{2}_{L^{2}}\nonumber\\
&=\varepsilon \|u^{-}\|^{\beta+1}_{L^{\beta+1}}+C(\varepsilon, m, \beta)\|u^{-}\|^{2}_{L^{2}},
\end{align}
where $\alpha=\frac{2(\beta+1-m)}{m(\beta-1)}$. \\
Similarly, for $\mu\in [2,\beta+1)$, we obtain
\begin{align}\label{levy}
\int_{Z}\psi(z)\nu(dz)\ \int_{D}{(u^{-})}^{\mu}(x,s)dxds&\leq C\|u^{-}\|^{\mu}_{L^{\mu}}\nonumber\\
&\leq C\|u^{-}\|^{\mu\alpha'}_{L^{2}}\|u^{-}\|^{\mu(1-\alpha')}_{L^{\beta+1}}\nonumber\\
&\leq  \varepsilon \|u^{-}\|_{L^{\beta+1}}^{\mu(1-\alpha')\frac{2}{2-\mu\alpha'}}+C(\varepsilon, \mu, \beta)\|u^{-}\|^{2}_{L^{2}}\nonumber\\
&=\varepsilon \|u^{-}\|^{\beta+1}_{L^{\beta+1}}+C(\varepsilon, \mu, \beta)\|u^{-}\|^{2}_{L^{2}},
\end{align}
where $\alpha'=\frac{2(\beta+1-\mu)}{\mu(\beta-1)}$. \\
Putting (\ref{beta}) and (\ref{levy}) into (\ref{u}), we obtain
\begin{align*}
\mathbb{E}\|\eta(u_{t})\|^{2}&\leq \int_{0}^{t}(2\varepsilon-2|a_{1}|)\mathbb{E}\|u^{-}_{s}\|^{\beta+1}_{L^{\beta+1}}ds+(2b_{2}+2a_{2}+C(\varepsilon, m, \beta)+C(\varepsilon, \mu, \beta))\int_{0}^{t}\mathbb{E}\|u_{s}^{-}\|^{2}_{L^{2}}ds
\end{align*}
Let $\varepsilon \in (0, |a_{1}|)$. Then
\begin{align*}
\mathbb{E}\|\eta(u_{t})\|^{2}\leq C \ \int_{0}^{t}\mathbb{E}\|\eta(u_{t})\|^{2}ds,
\end{align*}
From Gronwall's inequality, it follows that $\mathbb{E}\|\eta(u_{t})\|^{2}=0$. This implies that $\eta(u_{t})=u^{-}(x,t)=0$ a.s. for a.e. $x\in D$ and $t\in [0,T]$. The proof is complete. \ \ \ $\Box$ \\
\textbf{Remark 3.1.} The assumption $\mbox{(A1)}$ is weaker than the assumption  $\mbox{(H1)}$ in  \cite{Bao}. For example, if we consider the Allen-Cahn type equation, $f(u)=u-u^{3}$ doesn't satisfy $\mbox{(H1)}$, but $f$ satisfies $\mbox{(A1)}$. \\
\textbf{Remark 3.2.} Since $A=\sum_{i,j=1}^{d}\frac{\partial}{\partial x_{i}}[a_{ij}(x)\frac{\partial}{\partial x_{j}}]$ is more general than the Laplacian operator $\triangle$, Theorem 3.1 is the generalization of Theorem 3.1 in \cite{Lv}. \\
\textbf{Remark 3.3.}  If it is  assumed that $\beta\in (0,1)$, for the case $1+\beta\leq m<2$ and $1+\beta\leq q <2$, by the $L^{p}$ interpolation inequality and Young inequality, we can get the corresponding results.\\

\section{\textbf{Explosive solutions}}
In this section, we consider the unbounded solutions of the equation (\ref{rd}). \\
For the elliptic equation:
\begin{equation}\label{elliptic}
\left\{\begin{aligned}
&A\vartheta=-\lambda\vartheta,\ \mbox{in}\ D,\\
&\vartheta=0, \ \mbox{on}\ \partial D,
\end{aligned}\right.
\end{equation}
it is well known that all the eigenvalues of $-A$ are strictly positive, increasing and the eigenfunction $\phi$ corresponding to the smallest eigenvalue $\lambda_{1}$ does not change sign in the domain $D$ (see  p. 355, \cite{LE}). We can normalize it in such a way that
\begin{equation}\label{norm}
\phi(x)\geq 0, \ \ \ \int_{D}\phi(x)dx=1.
\end{equation}
\textbf{Theorem 4.1.}
Suppose the initial-boundary value problem (\ref{rd}) has a unique local solution and the conditions $\mbox{(A1)-(A4)}$ hold. In addition, we assume that  $\lambda_{1}> a_{2}$, $a_{1}>0$, and
\begin{equation*}
\int_{D}g(x)\phi(x)dx> \big(\frac{\lambda_{1}-a_{2}}{a_{1}}\big)^{\frac{1}{\beta-1}},
\end{equation*}
and if $\lambda_{1}<a_{2}$, we assume that $\int_{D}g(x)\phi(x)dx>0$, where $\lambda_{1}$ is the smallest eigenvalue of $-A$ and $\phi$ is the corresponding eigenfunction. Then, for any $p\geq 1$, there exists a constant $T_{p}>0$  such that
\begin{equation}\label{lp}
\lim_{t\rightarrow T^{-}_{p}}\mathbb{E}\|u_{t}\|_{L^{p}}=\lim_{t\rightarrow T^{-}_{p}}\mathbb{E}\big\{\int_{D}|u(x,t)|^{p}dx\big\}^{1/p}=\infty.
\end{equation}
That is, the solution explodes in mean $L^{p}$-norm sense. \\
\textbf{Proof.} By Theorem 3.1, Eq. (\ref{rd}) has a unique positive solution. We will prove the theorem by contradiction. We suppose (\ref{lp}) is false. Then there exists a global positive solution $u$ such that
\begin{equation*}
\sup_{0\leq t\leq T}\mathbb{E}\big\{\int_{D}|u(x,t)|^{p}dx\big\}^{1/p}<\infty,
\end{equation*}
for any $T>0$. Let $\phi$ be the the eigenfunction defined in (\ref{elliptic}). Define
\begin{equation}\label{hat}
\hat{u}(t):=\int_{D}u(x,t)\phi(x)dx\geq 0.
\end{equation}
By (\ref{norm}), $\phi$ can be regarded as the probability density function of a random variable $\xi$ in $D$, independent of $W_{t}$. The equality (\ref{hat})
can be interpreted as an expectation $\hat{u}(t)=\mathbb{E}_{\xi}\{u(\xi,t)\}$. From (\ref{rd}), (\ref{hat}) and the self-adjointness of $A$, it follows that
\begin{align}\label{multiply}
\hat{u}(t)&=(g,\phi)+\int_{0}^{t}\int_{D}[Au(x,s)]\phi(x)dxds+\int_{0}^{t}\int_{D}f(u,x,s)\phi(x)dxds\nonumber\\
&\quad+\int_{0}^{t}\int_{D}\sigma(u,\nabla u,x,s)\phi(x)dxdW(x,s)\nonumber\\
&\quad+\int_{0}^{t}\int_{Z}\int_{D}\varphi(u,x,z,s)\phi(x)dx\widetilde{N}(ds,dz)\nonumber\\
&=(g,\phi)-\lambda_{1}\int_{0}^{t}\int_{D}u(x,s)\phi(x)dxds+\int_{0}^{t}\int_{D}f(u,x,s)\phi(x)dxds\nonumber\\
&\quad+\int_{0}^{t}\int_{D}\sigma(u,\nabla u,x,s)\phi(x)dxdW(x,s)\nonumber\\
&\quad+\int_{0}^{t}\int_{Z}\int_{D}\varphi(u,x,z,s)\phi(x)dx\widetilde{N}(ds,dz).
\end{align}
Taking the expectation to both sides of (\ref{multiply}) and by Fubini's theorem, we have
\begin{align*}
\mathbb{E}\hat{u}(t)=(g,\phi)-\lambda_{1}\int_{0}^{t}\mathbb{E}\hat{u}(s)ds+\int_{0}^{t}\mathbb{E}\int_{D}f(u,x,s)\phi(x)dxds,
\end{align*}
or, in the differential form,
\begin{equation}\label{differential}
\left\{\begin{aligned}
&\frac{d\xi(t)}{dt}=-\lambda_{1}\xi(t)+\mathbb{E}\int_{D}f(u,x,s)dx\\
&\xi(0)=\xi_{0},
\end{aligned}\right.
\end{equation}
where $\xi(t)=\mathbb{E}\hat{u}(t)$, $\xi_{0}=(g,\phi)$. By $(\mbox{A1})$ and Jensen's inequality, we obtain
\begin{equation}\label{Jensen}
\left\{\begin{aligned}
&\frac{d\xi(t)}{dt}\geq-\lambda_{1}\xi(t)+a_{1}\xi^{\beta}(t)+a_{2}\xi(t),\\
&\xi(0)=\xi_{0},
\end{aligned}\right.
\end{equation}
If $\lambda_{1}\geq a_{2}$, for $\xi_{0}> \big(\frac{\lambda_{1}-a_{2}}{a_{1}}\big)^{\frac{1}{\beta-1}}$, we can show that $\xi(\cdot)$ is strictly increasing. It follows from (\ref{Jensen}) that
\begin{align}\label{contra}
T&\leq \int_{\xi_{0}}^{\xi(T)}\frac{ds}{a_{1}s^{\beta}-(\lambda_{1}-a_{2})s}\leq \int_{\big(\frac{\lambda_{1}-a_{2}}{a_{1}}\big)^{\frac{1}{\beta-1}}}^{\infty}\frac{ds}{a_{1}s^{\beta}-(\lambda_{1}-a_{2})s}< \infty.
\end{align}
If $\lambda_{1}< a_{2}$, for $\xi_{0}> 0$, we can show that $\xi(\cdot)$ is strictly increasing. We have
\begin{align}\label{contradict}
T=\int_{0}^{T}dt&\leq \int_{0}^{T}\frac{d\xi(t)}{a_{1}\xi^{\beta}(t)}=\int_{\xi_{0}}^{\xi(T)}\frac{ds}{a_{1}s^{\beta}}< \infty.
\end{align}
Since $T$ is arbitrary, either (\ref{contra}) or (\ref{contradict}) results in a contradiction. Therefore, for $\lambda_{1}\geq a_{2}$, $\xi_{0}> \big(\frac{\lambda_{1}-a_{2}}{a_{1}}\big)^{\frac{1}{\beta-1}}$,  $\xi(t)$ must blow up at a time $T_{p}\leq \int_{\xi_{0}}^{\xi(T)}\frac{ds}{a_{1}s^{\beta}-(\lambda_{1}-a_{2})s}$. For $\lambda_{1}< a_{2}$, $\xi(t)$ must blow up at a time $T_{p}\leq \int_{\xi_{0}}^{\xi(T)}\frac{ds}{a_{1}s^{\beta}}$. \\
Since $\phi$ is bounded and continuous on $\overline{D}$, by H$\ddot{\mbox{o}}$lder's inequality, we have
\begin{align}
\xi(t)\leq \big(\int_{D}|\phi(x)|^{q}dx \big)^{1/q}\big(\mathbb{E}\int_{D}|u(x,t)|^{p}dx \big)^{1/p},
\end{align}
where $q=p/(p-1)$, $p\geq 1$. So the positive solution explodes at some time $T'\leq T_{e}$ in the mean $L^{p}$-norm for each $p\geq 1$. The proof is complete. \ \ \ $\Box$ \\
\textbf{Example 4.1.}
Consider the following problem in a spherical domain $D=B(R)$ in $\mathbb{R}^{3}$:
\begin{equation}\label{example}
\left\{\begin{aligned}
&\frac{\partial u}{\partial t}=\triangle u+u^{\frac{8}{3}}-u+\gamma_{0}(u^{3}+|\nabla u|^{2})^{1/2}\partial_{t}W(x,t)+c_{0}\int_{0}^{\infty}zu^{3}\partial_{t}\widetilde{N}(t,dz),\
t>0,x\in D,\\
&u(x,0)=a_{0}e^{-\alpha |x|},\ x\in D,\\
&u(x,t)|_{|x|=R}=0, \ t>0,
\end{aligned}\right.
\end{equation}
where $\widetilde{N}(dt,dz)=N(dt,dz)-dt\nu(dz)$ is a compensated Poisson measure corresponding to the Poisson random measure $N(dt,dz)$, $W(x,t)$ is a continuous Wiener random field with the covariance function
\begin{align*}
q(x,y)=b_{0}\exp\{-\rho(x\cdot y)\}, \ \ x,y\in \mathbb{R}^{3}.
\end{align*}
The constants $\sigma_{0}$, $c_{0}$, $a_{0}$, $\alpha$ are strictly positive and $x\cdot y=\sum_{i=1}^{3}x_{i}y_{i}$. Here $A=\triangle$, $f=u^{\frac{8}{3}}-u$, $\sigma=\gamma_{0}(u^{3}+|\nabla u|^{2})^{1/2}$, $\varphi=c_{0}zu^{3}$, $Z=(0,\infty)$.
It is obvious that conditions $\mbox{(A1)}$ and  $\mbox{(A4)}$ are  satisfied. If $\frac{1}{2}b_{0}\gamma^{2}_{0}<1$,
since
\begin{align*}
\frac{1}{2}b_{0}\gamma^{2}_{0}\exp\{-\rho|x|^{2}\}(u^{3}+|\xi|^{2})-|\xi|^{2}\leq (\frac{1}{2}b_{0}\gamma^{2}_{0}-1)|\xi|^{2}+\frac{1}{2}b^{2}\gamma^{2}_{0}u^{3}.
\end{align*}
 then condition $\mbox{(A2)}$ is satisfied.
 If $\int_{0}^{\infty}z^{2}\nu(dz)<\infty$, take $\mu=6$ and $\psi(z)=c^{2}_{0}z^{2}$, then the condition $\mbox{(A3)}$ is satisfied.
From Theorem 3.1, it follows that the solution of Eq. (\ref{example}) is positive. The smallest eigenvalue of the elliptic equation (\ref{elliptic}) is $\lambda_{1}=(\frac{\pi}{R})^{2}$ and the corresponding normalized eigenfunction is $\phi(x)=\frac{1}{4R^{2}|x|}\sin\frac{\pi|x|}{R}$, $0<|x|<R$.
If $a_{0}$ is sufficiently large, then we have
\begin{align}\label{a0}
\int_{D}g(x)\phi(x)dx=\int_{0}^{R}\frac{a_{0}e^{-\alpha r}}{4R^{2}r}\sin \frac{\pi r}{R}dr >\frac{a_{0}}{4R^{3}}\int_{0}^{R}e^{-\alpha r}\sin\frac{\pi r}{R}dr>(\frac{\pi^{2}}{R^{2}}+1)^{\frac{3}{5}},
\end{align}
Therefore, by Theorem 4.1, the solutions to the Eq. (\ref{example}) will blow up in finite time in mean $L^{p}$-norm for any $p\geq 1$.
Note that Theorem 3.1 in \cite{Bao} is not suitable for Eq. (\ref{example}).

To discuss the noise-induced explosion, we consider the following stochastic reaction-diffusion equation:
\begin{equation}\label{noise-induced}
\left\{\begin{aligned}
&\frac{\partial u}{\partial t}=Au+f(u,x,t)+\sigma(u,x,t)\partial_{t}W(x,t)\\
&\quad\quad\quad+\int_{Z}\varphi(u,x,z,t)\partial_{t}\widetilde{N}(t,dz),\
t>0,x\in D\\
&u(x,0)=g(x),\ x\in D,\\
&u(x,t)=0, \ t>0,x\in \partial D,
\end{aligned}\right.
\end{equation}
which is a special case of Eq. (\ref{rd}), where $\sigma$ is independent of $\nabla u$. We assume that the noise terms satisfy the following conditions: \\
$(\mbox{A1}')$ The correlation function $q(x,y)$ is continuous and positive for $x,y\in \overline{D}$ such that
\begin{align*}
\int_{D}\int_{D}q(x,y)v(x)v(y)dxdy\geq \kappa \big (\int_{D}v(x)dx\big)^{2}
\end{align*}
for any non-negative $v\in H$ and some constant $\kappa>0$. \\
$(\mbox{A2}')$  The function $f(u,x,t)$ is continuous on $\mathbb{R}\times \overline{D} \times [0,\infty)$ such that $f(u,x,t)\geq 0$ for $u\geq 0$ and $x\in \overline{D}$, $t\in [0,\infty)$. \\
$(\mbox{A3}')$ There exist continuous functions $\sigma_{0}$, $G$ such that they are both positive, convex and satisfy
\begin{align*}
\sigma(u,x,t)\geq \sigma_{0}(u), \ \ \sigma^{2}_{0}(u)\geq G(u^{2}),
\end{align*}
for $x\in \overline{D}$, $u\geq 0$, $t\in [0,\infty)$. \\
$(\mbox{A4}')$ There exist continuous functions $\varphi_{0}$, $K$ such that they are both positive, convex  and satisfy
\begin{align*}
&\int_{Z}\big(\int_{D}\varphi(u,x,z,t)\phi(x)dx\big)^{2}\nu(dz)\geq \int_{Z}\big(\int_{D}\varphi_{0}(u,z)\phi(x)dx\big)^{2}\nu(dz), \\ &\int_{Z}\varphi^{2}_{0}(u,z)\nu(dz)\geq K(u^{2}),
\end{align*}
for $x\in \overline{D}$, $z\in Z$, $u\geq 0$, $t\in [0,\infty)$. \\
$(\mbox{A5}')$ There exists a constant $M>0$ such that
$\kappa G(u)+K(u)> 2\lambda_{1}u$ for $u>M$, and
\begin{align*}
\int_{M}^{\infty}\frac{du}{\kappa G(u)+K(u)-2\lambda_{1}u}< \infty.
\end{align*}
$(\mbox{A6}')$ The initial datum satisfies the following
\begin{align*}
(g,\phi)=\int_{D}g(x)\phi(x)dx>M.
\end{align*}
\textbf{Theorem 4.2.} Assume that the initial-boundary value problem (\ref{rd}) has a unique positive local solution and the conditions $(\mbox{A1}')-(\mbox{A6}')$ hold. Then for each $p\geq 2$, there exists a constant $T_{p}$ such that
\begin{equation}\label{noise}
\lim_{t\rightarrow T^{-}_{p}}\mathbb{E}\|u_{t}\|_{L^{p}}=\lim_{t\rightarrow T^{-}_{p}}\mathbb{E}\big\{\int_{D}|u(x,t)|^{p}dx\big\}^{1/p}=\infty.
\end{equation}
That is, the solution explodes in mean $L^{p}$-norm sense. \\
\textbf{Proof.}  We assume the conclusion is false. Then there exists the solution $u$ and for some $p\geq 2$, $\mathbb{E}\|u_{t}\|^{p}<\infty$, $t\in [0,T]$, for any $T>0$.
Let $\hat{u}(t)=(\phi,u_{t})$ be defined as before. By (\ref{noise-induced}),
\begin{align}\label{gradient}
\hat{u}(t)&=(g,\phi)+\int_{0}^{t}\int_{D}[Au(x,s)]\phi(x)dxds+\int_{0}^{t}\int_{D}f(u,x,s)\phi(x)dxds\nonumber\\
&\quad+\int_{0}^{t}\int_{D}\sigma(u,x,s)\phi(x)dW(x,s)dx\nonumber\\
&\quad+\int_{0}^{t}\int_{Z}\int_{D}\varphi(u,x,z,s)\phi(x)dx\widetilde{N}(ds,dz)\nonumber\\
&=(g,\phi)-\lambda_{1}\int_{0}^{t}\int_{D}u(x,s)\phi(x)dxds+\int_{0}^{t}\int_{D}f(u,x,s)\phi(x)dxds\nonumber\\
&\quad+\int_{0}^{t}\int_{D}\sigma(u,x,s)\phi(x)dW(x,s)dx\nonumber\\
&\quad+\int_{0}^{t}\int_{Z}\int_{D}\varphi(u,x,z,s)\phi(x)dx\widetilde{N}(ds,dz).
\end{align}
By (\ref{gradient}), we apply the It$\hat{\mbox{o}}$'s formula to $\hat{u}^{2}(t)$ to get
\begin{align}\label{square}
\hat{u}^{2}(t)&=(g,\phi)^{2}-2\lambda_{1}\int_{0}^{t}\hat{u}^{2}(s)ds+2\int_{0}^{t}\int_{D}\hat{u}(s)f(u,x,s)\phi(x)dxds\nonumber\\
&\quad+2\int_{0}^{t}\int_{D}\hat{u}(s)\sigma(u,x,s)\phi(x)dxdW(x,s)\nonumber\\
&\quad+\int_{0}^{t}\int_{D}\int_{D}q(x,y)\phi(x)\phi(y)\sigma(u,x,s)\sigma(u,y,s)dxdyds\nonumber\\
&\quad+\int_{0}^{t}\int_{Z}\big[\big(\hat{u}(s)+\int_{D}\varphi(u,x,z,s)\phi(x)dx\big)^{2}-\hat{u}^{2}(s)\big]\widetilde{N}(ds,dz)\nonumber\\
&\quad+\int_{0}^{t}\int_{Z}\big[\big(\hat{u}(s)+\int_{D}\varphi(u,x,z,s)\phi(x)dx\big)^{2}-\hat{u}^{2}(s)\nonumber\\
&\quad-2\hat{u}(s)\int_{D}\varphi(u,x,z,s)\phi(x)dx\big]\nu(dz)ds.
\end{align}
Let $\eta(t)=\mathbb{E}\hat{u}^{2}(t)$.  Taking expectations of both sides of (\ref{square}), we obtain
\begin{align*}
\eta(t)&=(g,\phi)^{2}-2\lambda_{1}\int_{0}^{t}\eta(s)ds+2\mathbb{E}\int_{0}^{t}\int_{D}\hat{u}(s)f(u,x,s)\phi(x)dxds\nonumber\\
&\quad+\mathbb{E}\int_{0}^{t}\int_{D}\int_{D}q(x,y)\phi(x)\phi(y)\sigma(u,x,s)\sigma(u,y,s)dxdyds\nonumber\\
&\quad+\mathbb{E}\int_{0}^{t}\int_{Z}\big(\int_{D}\varphi(u,x,z,s)\phi(x)dx\big)^{2}\nu(dz)ds,
\end{align*}
or, in the differential form,
\begin{equation}\label{form}
\left\{\begin{aligned}
&\frac{d\eta(t)}{dt}=-2\lambda_{1}\eta(t)+2\mathbb{E}\hat{u}(t)\int_{D}f(u,x,t)\phi(x)dx\\
&\quad\quad\quad+\mathbb{E}\int_{D}\int_{D}q(x,y)\phi(x)\phi(y)\sigma(u,x,t)\sigma(u,y,t)dxdy\\
&\quad\quad\quad+\mathbb{E}\int_{Z}\big(\int_{D}\varphi(u,x,z,t)\phi(x)dx\big)^{2}\nu(dz), \\
&\eta(0)=\eta_{0}=(g,\phi)^{2}.
\end{aligned}\right.
\end{equation}
By $(\mbox{A2}')$, we have $\mathbb{E}\hat{u}(t)\int_{D}f(u,x,t)\phi(x)dx\geq 0$. By $(\mbox{A1}')$, $(\mbox{A3}')$, Jensen's inequality, we have
\begin{align}\label{Schwarz}
&\int_{D}\int_{D}q(x,y)\phi(x)\phi(y)\sigma(u,x,t)\sigma(u,y,t)dxdy\nonumber\\
&\geq \kappa\big(\int_{D}\phi(x)\sigma(u,x,t)dx\big)^{2}\nonumber\\
&\geq \kappa\big(\int_{D}\phi(x)\sigma_{0}(u)dx\big)^{2}\nonumber\\
&\geq \kappa\sigma^{2}_{0}(\hat{u}(t))\nonumber\\
&\geq \kappa G(\hat{u}^{2}(t)).
\end{align}
By $(\mbox{A4}')$,  Jensen's inequality, we get
\begin{align}\label{dz}
&\int_{Z}(\int_{D}\varphi(u,x,z,t)\phi(x)dx)^{2}\nu(dz)\nonumber\\
&\geq\int_{Z}\big(\int_{D}\varphi_{0}(u,z)\phi(x)dx\big)^{2}\nu(dz)\nonumber\\
&\geq\int_{Z}\varphi^{2}_{0}(\hat{u}(t),z)\nu(dz)\nonumber\\
&\geq K(\hat{u}^{2}(t)).
\end{align}
From (\ref{Schwarz}), (\ref{dz}), (\ref{form}) and Jensen's inequality,  it follows that
\begin{align}\label{deduce}
\frac{d\eta(t)}{dt}&\geq -2\lambda_{1}\eta(t)+\kappa\mathbb{E}G(\hat{u}^{2}(t))+\mathbb{E}K(\hat{u}^{2}(t))\nonumber\\
&\geq -2\lambda_{1}\eta(t)+\kappa G(\eta(t))+K(\eta(t)).
\end{align}
Similar to (\ref{contra}), for $\eta_{0}>M$, we obtain
\begin{align*}
T\leq \int_{\eta_{0}}^{\eta(T)}\frac{du}{\kappa G(\eta(t))+K(\eta(t))-2\lambda_{1}u}\leq\int_{M}^{\infty}\frac{du}{\kappa G(\eta(t))+K(\eta(t))-2\lambda_{1}u}<\infty.
\end{align*}
Since $T$ is arbitrary, this results in a contradiction. Therefore, the mean square $\eta(t)=\mathbb{E}\hat{u}^{2}(t)$ must blow up at some finite time $T_{\ast}>0$. Applying the H$\ddot{\mbox{o}}$lder inequality, we see that (\ref{noise}) holds for each $p\geq 2$. \ \ $\Box$ \\
 \textbf{Remark 4.1.} In \cite{Ch},  \cite{CL}, the correlation function $q(x,y)$ is assumed to satisfy the inequality
 \begin{equation}\label{remark}
 \int_{D}\int_{D}q(x,y)v(x)v(y)dxdy\geq q_{1}\int_{D}v^{2}(x)dx
 \end{equation}
 for any positive $v \in H$ and for some $q_{1}>0$.

 In fact, this assumption is not suitable. If the domain $D$ is bounded and $q\in (0,1]$, by the Cauchy-Schwarz inequality, we have
\begin{equation}\label{correlation}
\int_{D}\int_{D}q(x,y)v(x)v(y)dxdy\leq \big(\int_{D}v(x)dx\big)^{2}\leq \mu(D)\int_{D}v^{2}(x)dx,
\end{equation}
where $\mu(D)$ is the volume of $D$. By (\ref{remark}), (\ref{correlation}), we have $\mu(D)\geq q_{1}$. If the bounded domain $D$ is small enough, then we get a contradiction. \\
 \textbf{Remark 4.2.}  We consider the problem (\ref{noise-induced}) with a Levy-type noise, and the coefficient operator $A$ is more general than the Laplacian, it is easy to see that   Theorem 4.2 is the generalization of Theorem 4.3 in \cite{Lv}. \\
\textbf{Example 4.2.}
Consider the following problem in a spherical domain $D=B(R)$ in $\mathbb{R}^{3}$:
\begin{equation}\label{noise explosion}
\left\{\begin{aligned}
&\frac{\partial u}{\partial t}=\triangle u+u^{1+\alpha}+\mu u^{4}\partial_{t}W(x,t)+c_{0}\int_{0}^{\infty}zu^{6}\partial_{t}\widetilde{N}(t,dz),\
t>0,x\in D,\\
&u(x,0)=a_{0}e^{-\beta |x|},\ x\in D,\\
&u(x,t)|_{|x|=R}=0, \ t>0,
\end{aligned}\right.
\end{equation}
where $\widetilde{N}(dt,dz)=N(dt,dz)-dt\nu(dz)$ is a compensated Poisson measure corresponding to the Poisson random measure $N(dt,dz)$, $W(x,t)$ is a continuous Wiener random field with the covariance function
\begin{align*}
q(x,y)=b_{0}\exp\{-\rho(x\cdot y)\}, \ \ x,y\in \mathbb{R}^{3}.
\end{align*}
The constants $\mu$, $c_{0}$, $a_{0}$, $\alpha$, $\beta$ are strictly positive and $x\cdot y=\sum_{i=1}^{3}x_{i}y_{i}$. Here $A=\triangle$, $f=u^{1+\alpha}$, $\sigma=\mu u^{4}$, $\varphi=c_{0}zu^{6}$, $Z=(0,\infty)$.\\
For $x,y\in B(R)$, we have
\begin{align*}
q(x,y)\geq \kappa=b_{0}\exp\{-\rho R^{2}\}.
\end{align*}
Then for any non-negative $v\in H$,
\begin{align*}
\int_{D}\int_{D}q(x,y)v(x)v(y)dxdy\geq \kappa \big (\int_{D}v(x)dx\big)^{2}.
\end{align*}
The condition $(\mbox{A1}')$  holds. It is obvious that $f=u^{1+\alpha}> 0$ for $u>0$, the condition $(\mbox{A2}')$ holds.
Let $G(u)=\mu^{2} u^{4}$, $\sigma_{0}(u)=\mu u^{4}$. Then $\sigma(u,x,t)=\mu u^{4}=\sigma_{0}(u)$ and $\sigma^{2}_{0}(u)=G(u^{2})$. $\sigma_{0}$ and $G$ are both continuous, positive and convex. The condition $(\mbox{A3}')$ is satisfied. Let $\varphi_{0}(u,z)=c_{0}zu^{6}$. Assume $\int_{0}^{\infty}z^{2}\nu(dz)<\infty$,  let $K(u)=\big(c_{0}\int_{0}^{\infty}z^{2}\nu(dz)\big)u^{6}$. Then $\varphi(u,x,z,t)=c_{0}zu^{6}=\varphi_{0}(u,z)$. $\varphi_{0}$ and $K$ are both positive and convex. The condition $(\mbox{A4}')$ holds. The smallest eigenvalue of the elliptic equation (\ref{elliptic}) is $\lambda_{1}=(\frac{\pi}{R})^{2}$. If $b_{0}$ or $c_{0}$ is large enough, it is easy to see that $b_{0}\exp\{-\rho R^{2}\}\mu^{2}u^{4}+\big(c_{0}\int_{0}^{\infty}z^{2}\nu(dz)\big)u^{6}>2(\frac{\pi}{R})^{2}u$ for $u> M$, and
\begin{align*}
\int_{M}^{\infty}\frac{du}{b_{0}\exp\{-\rho R^{2}\}\mu^{2}u^{4}+\big(c_{0}\int_{0}^{\infty}z^{2}\nu(dz)\big)u^{6}-2(\frac{\pi}{R})^{2}u}< \infty.
\end{align*}
The condition  $(\mbox{A5}')$  is satisfied. If $a_{0}$ is sufficiently large, simliar to (\ref{a0}),we have $\int_{D}g(x)\phi(x)dx>M$.
By Theorem (4.2), the solution of (\ref{noise explosion}) will blow up in finite time in $L^{p}$-norm for any $p\geq 2$.

Now we consider the explosive solution problem for (\ref{rd}), when the domain $D$ is unbounded, for example, $D=\mathbb{R}^{d}$. Let $B(R)=\{x\in \mathbb{R}^{d}: |x|< R\}$. \\
 \textbf{Theorem 4.3.}  Assume that the initial-boundary value problem (\ref{rd}) has a unique local solution and the conditions $\mbox{(A1)-(A3)}$ hold, where $D=\mathbb{R}^{d}$. Then for any $R>0$,  there exists a constant $T_{p}(R)>0$  such that
\begin{align*}
\lim_{t\rightarrow T_{p}(R)-}\mathbb{E}\big\{\int_{B(R)}|u(t,x)|^{p}\big\}^{1/p}=\infty,
\end{align*}
provided that the conditions of Theorem 4.1 holds for $p\geq 1$ or the conditions of Theorem 4.2 holds for $p\geq 2$, where $D=\mathbb{R}^{d}$. \\
\textbf{Proof.} The proof follows the spirit of the one for Theorem 3.3 in \cite{CL}. For the sake of completeness, we present it. We only consider the case under the conditions of Theorem 4.1, since the proof under the conditions of Theorem 4.2 is similar.

By restricting the solution $u$ to $\overline{B(R)}$, let $\hat{u}(t)=\int_{B(R)}u(x,t)\phi(x)dx\geq 0$ as defined by (\ref{hat}). Since $u\geq 0$ on the boundary $\partial B(R)$, by Green's identity,
\begin{align}\label{normal}
(Au_{t},\phi)=-\lambda_{1}\phi+\int_{\partial B(R)}u(x,t)\big(-\frac{\partial \phi(x)}{\partial \nu}\big)dS.
\end{align}
 Denote $n=(n_{1},n_{2},\ldots,n_{d})$ as the unit outward normal vector to the boundary $\partial B(R)$. Since there exists a constant $c>0$ such that $\sum_{i,j=1}^{d}a_{ij}(x)\xi_{i}\xi_{j}\geq c|\xi|^{2}$ for all $x\in \overline{D}$ and $\xi=(\xi_{1},\ldots, \xi_{d})\in \mathbb{R}^{d}$. We have $\nu\cdot n=\sum_{i,j}^{d}a_{ij}n_{i}n_{j}\geq 0$. This implies that the conormal $\nu(x)$ is an exterior direction field. Since $\phi>0$ in $B(R)$ and $\phi=0$ on $\partial B(R)$, we have
 \begin{align}\label{direction}
\frac{\partial \phi(x)}{\partial \nu}\leq 0.
 \end{align}
 Putting (\ref{normal}) into (\ref{multiply}), by (\ref{direction}), we obtain the inequality (\ref{Jensen}). The rest of proof can be completed as in Theorem 4.1. \ \  $\Box$ \\
\textbf{Remark 4.3.} In \cite{Bao}, when the domain $D=\mathbb{R}^{d}$, it seems impossible to consider the existence of the position solution of initial-boundary value problem (\ref{rd}) and the explosionn of the position solution. The reason is that the proof of Theorem 2.1 in \cite{Bao} relies on the fact that the volume $V(D)$ of domain $D$ is bounded. In \cite{Bao}, the proofs of Theorem 3.1 and Theorem 4.1 are both rely on Theorem 2.1. So for $D=R^{d}$, they are not valid.
\section{\textbf{Global solutions for a stochastic Allen-Cahn  equation driven by a L$\acute{\mbox{e}}$vy type noise}}
In this section, we consider the following stochastic Allen-Cahn equation driven by a L$\acute{\mbox{e}}$vy type noise,
\begin{equation}\label{AC}
\left\{\begin{aligned}
&du=(\triangle u+u(1-u^{2}))dt+bu^{m}dW_{t}+cu^{n}\int_{Z}z\widetilde{N}(dt,dz), \ t>0, \ x\in D,\\
&u(x,0)=h(x),\ x\in D,\\
&u(x,t)=0, \ t>0,x\in \partial D,
\end{aligned}\right.
\end{equation}
where $1<m<2$, $1<n<2$, $b,c\in \mathbb{R}$, $D\subset \mathbb{R}^{3}$, $Z=(0,\infty)$.

Let $V$ be a real separable Hilbert space. We first consider the more general equation
\begin{equation}\label{general}
\left\{\begin{aligned}
&du_{t}=(Au_{t}+F_{t}(u_{t})dt+\Sigma_{t}(u_{t})dW_{t}+\int_{Z}\Gamma_{t}(u_{t},z)\widetilde{N}(dt,dz),\ t\geq 0,\\
&u_{0}=h(x),
\end{aligned}\right.
\end{equation}
where the coefficients $A$, $F_{t}$, $\Sigma_{t}$ and $\Gamma_{t}$ are assumed to be non-random or deterministic. $W(x,t)$ is a Wiener random field, $(Z,\mathcal{B}(Z))$ is a measurable space. $\widetilde{N}(dt,dz)$ is the compensated Poisson measure. Here we say that an $\mathcal{F}_{t}$-adapted $V$-valued process $u_{t}$ is a \emph{strong solution}, or a \emph{variational solution}, of the equation (\ref{general}) if $u\in L^{2}([0,T]; V)$, and for any $\varphi\in V$, the following equation
\begin{align*}
(u_{t},\varphi)&=(h,\varphi)+\int_{0}^{t}(Au_{s}, \varphi)ds+\int_{0}^{t}(F_{s}(u_{s}), \varphi)ds+\int_{0}^{t}(\varphi, \Sigma_{s}(u_{s})dW_{s})\\
&\quad+\int_{0}^{t}(\int_{Z}\Gamma_{s}(u_{s},z)\widetilde{N}(ds,dz), \varphi)
\end{align*}
holds for each $t\in [0,T]$ a.s.

Denote $L_{1}(V)$ the space of nuclear (trace class) operators on $V$. Let $U\subset V$ be a open set and let $U\times [0,T]=U_{T}$. Here a functional $\Phi: U_{T}\rightarrow \mathbb{R}$ is said to be a \emph{strong It$\hat{\mbox{o}}$ functional} if it satisfies the following (see \cite[pp. 226]{C}): \\
(1) $\Phi: U_{T}\rightarrow \mathbb{R}$ is locally bounded and continuous such that its first two partial derivatives $\partial_{t}\Phi(v,t)$, $\Phi'(v,t)$ and $\Phi''(v,t)$ exist for each $(v,t)\in U_{T}$.\\
(2) The derivatives $\partial _{t}\Phi$ and $\Phi'\in V$ are locally bounded and continuous in $U_{T}$.\\
(3) For any $\Gamma\in \mathcal{L}_{1}(V)$, the map: $(v,t)\rightarrow Tr[\Phi''(v,t)\Gamma]$ is locally bounded and continuous in $(v,t)\in U_{T}$.\\
(4) $\Phi'(\cdot,t): U\rightarrow V$ is such that $(\Phi'(\cdot,t), v)$ is continuous in $t\in [0,T]$ for any $v\in V$ and
\begin{align*}
\|\Phi'(v,t)\|\leq \kappa(1+\|v\|), \  (v,t)\in U\times [0,T],
\end{align*}
for some $\kappa>0$.

Let $U\subset V$ be a neighborhood of the origin. By a similar statement to that in \cite[pp. 228]{C}, we present the definition of Lyapunov functional.
Define the operator $\mathcal{L}_{t}$ as follows:
\begin{align*}
\mathcal{L}_{t}\Phi(v,t)&=\frac{\partial}{\partial t}\Phi(v,t)+\frac{1}{2}Tr[\Phi''(v,t)\Sigma_{t}(v)Q\Sigma^{\ast}_{t}(v)]+(Av, \Phi'(v,t))\\
&\quad+(F_{t}(v), \Phi'(v,t))+\int_{Z}[\Phi(v+\Gamma_{t}(v,z),t)-\Phi(v,t)-\big(\Gamma_{t}(v,z), \Phi'(v,t)\big)]\nu(dz),
\end{align*}
where $Q$ is the covariance operator.

A strong It$\hat{\mbox{o}}$ functional $\Phi: U\times \mathbb{R}^{+}\rightarrow \mathbb{R}$ is said to be a
\emph{Lyapunov functional }for the equation (\ref{general}), if \\
(1) $\Phi(0,t)=0$ for all $t\geq 0$, and, for any $\varepsilon>0$, there is a $\delta>0$ such that
\begin{align*}
\inf_{t\geq 0, \|h\|\geq \varepsilon}\Phi(h,t)\geq \delta, \ \mbox{and}
\end{align*}
(2) for any $t\geq 0$ and $v\in U$,
\begin{align*}
\mathcal{L}_{t}\Phi(v,t)\leq 0.
\end{align*}

Let $u^{h}_{t}$ be a strong solution of the equation (\ref{AC}) with $u^{h}_{0}=h$. \\
\textbf{Definition 5.1.} The solution $u^{h}_{t}$ is said to be nonexplosive if
\begin{equation*}
\lim_{r\rightarrow \infty} P\{\sup_{0\leq t\leq T}\|u^{h}_{t}\|\geq r\}=0,
\end{equation*}
for any $T>0$. If the above holds for $T=\infty$, the solution is said to be ultimately bounded. \\
\textbf{Lemma 5.1.} Let $\Phi:U\times \mathbb{R}^{+}\rightarrow \mathbb{R}^{+}$ be a Lyapunov functional and let $u^{h}_{t}$ denote the strong solution of (\ref{general}). For $r>0$, let $B_{r}=\{h\in V: \|h\|<r\}$ such that $B_{r}\subset U$. Define
\begin{align*}
\tau=\inf\{t>0: u^{h}_{t}\in B^{c}_{r}, \ h\in B_{r}\},
\end{align*}
with $ B^{c}_{r}=V\backslash B_{r}$. We put $\tau=T$ if the set is empty. Then the process $\phi_{t}=\Phi(u^{h}_{t\wedge \tau}, t\wedge \tau)$ is a local $\mathcal{F}_{t}$-supermartingale and the following Chebyshev inequality holds
\begin{align*}
P\{\sup_{0\leq t \leq T}\|u^{h}_{t}\|\geq r\}\leq \frac{\Phi(h,0)}{\Phi_{r}},
\end{align*}
where
\begin{align*}
\Phi_{r}=\inf_{0\leq t \leq T, h\in U\cap B^{c}_{r}} \Phi(h,t).
\end{align*}
\textbf{Proof.} From It$\hat{\mbox{o}}$'s formula, it follows that
\begin{align*}
\Phi(u^{h}_{t\wedge \tau})&=\Phi(h,0)+\int_{0}^{t\wedge \tau}\mathcal{L}_{s}\Phi(u^{h}_{s},s)ds+\int_{0}^{t\wedge \tau}(\Phi'(u^{h}_{s},s),\Sigma_{s}(u^{h}_{s}))dW_{s}\\
&\quad+\int_{0}^{t\wedge \tau}\int_{Z}(\Phi(v+\Gamma_{s}(v,z),s)-\Phi(v,s))\widetilde{N}(ds,dz)\\
&\leq \Phi(h,0)+\int_{0}^{t\wedge \tau}(\Phi'(u^{h}_{s},s),\Sigma_{s}(u^{h}_{s}))dW_{s}\\
&\quad+\int_{0}^{t\wedge \tau}\int_{Z}(\Phi(v+\Gamma_{s}(v,z),s)-\Phi(v,s))\widetilde{N}(ds,dz),
\end{align*}
therefore, $\phi_{t}=\Phi({u^{h}_{t\wedge \tau}}, t\wedge \tau)$ is a local supermartingale and
\begin{align*}
\mathbb{E}\phi_{t}\leq \mathbb{E}\phi_{0}=\Phi(h,0).
\end{align*}
By definition,
\begin{align*}
\mathbb{E}\phi_{T}&=\mathbb{E}\Phi(u^{h}_{T\wedge \tau}, T\wedge \tau)\\
&\geq \mathbb{E}\{\Phi(u^{h}_{\tau}; \tau\leq T)\}\\
&\geq \inf_{0\leq t\leq T, \|h\|=r}\Phi(h,t)P\{\tau\leq T\}\\
&\geq \Phi_{r}P\{\sup_{0\leq t\leq T}\|u^{h}_{t}\|\geq r\},
\end{align*}
the proof is complete.  \ \ $\Box$ \\
\textbf{Theorem 5.1.} If there exists an It$\hat{\mbox{o}}$ functional $\Psi: V\times \mathbb{R}^{+}\rightarrow \mathbb{R}^{+}$ and a constant $\alpha>0$ such that
\begin{align*}
\mathcal{L}_{t}\Psi\leq \alpha \Psi(v,t)\ \  \mbox{for any} \ v\in V,
\end{align*}
and the infimum  $\inf_{t\geq 0, \|h\|\geq r}\Psi(h,t)=\Psi_{r}$ exists such that $\lim_{r\rightarrow \infty}\Psi_{r}=\infty$, then the solution $u^{h}_{t}$ does not explode in finite time. \\
\textbf{Proof.} Let $\Phi(v,t)=e^{-\alpha t}\Psi(v,t)$. We have
\begin{align*}
\mathcal{L}_{t}\Phi(v,t)&=\frac{\partial}{\partial t}\Phi(v,t)+\frac{1}{2}\mbox{Tr}[\Phi''(v,t)\Sigma_{t}(v)Q\Sigma^{\ast}_{t}(v)]+(Av, \Phi'(v,t))\\
&\quad+(F_{t}(v), \Phi'(v,t))+\int_{Z}[\Phi(v+\Gamma_{t}(v,z),t)-\Phi(v,t)-\Phi'(v,t)\Gamma_{t}(v,z)]\nu(dz)\\
&=-\alpha e^{-\alpha t}\Psi(v,t)+e^{-\alpha t}\frac{\partial}{\partial t}\Psi(v,t)+e^{-\alpha t}\big(\frac{1}{2}Tr[\Psi''(v,t)\Sigma_{t}(v)Q\Sigma^{\ast}_{t}(v)]+(Av, \Psi'(v,t))\\
&\quad+(F_{t}(v), \Psi'(v,t))+\int_{Z}[\Psi(v+\Gamma_{t}(v,z),t)-\Psi(v,t)-\Psi'(v,t)\Gamma_{t}(v,z)]\nu(dz)\big)\\
&=-\alpha e^{-\alpha t}\Psi(v,t)+e^{-\alpha t}\mathcal{L}_{t}\Psi\leq 0.
\end{align*}
Therefore $\Phi$ is a Lyapunov functional. By Lemma 5.1, we have
\begin{align*}
P\{\sup_{0\leq t \leq T}\|u^{h}_{t}\|>r\}\leq \frac{\Phi(h,0)}{\Phi_{r}}=\frac{\Psi(h,0)}{\Psi_{r}}\rightarrow 0
\end{align*}
as $r\rightarrow \infty$, for any $T>0$.  \ \ $\Box$ \\
\textbf{Theorem 5.2.} Let $1<m<2$, $1<n<2$ and the initial datum $h(x)$ on $\overline{D}$ is positive and continuous. Suppose that $\int_{0}^{\infty}z^{2}\nu(dz)< \infty$ and there is $q_{0}>0$ such that $\sup_{x,y\in D}|q(x,y)|\leq q_{0}$. Then the equation (\ref{AC}) has a global strong solution. \\
\textbf{Proof.} In view of  the  proof of Theorem 3-6.5 in \cite[pp. 86]{C} and the proof of Theorem 3.2 in \cite{RZ}, we can show that the equation (\ref{AC})  has a local strong solution. By Theorem 3.1, the solution is positive. Define $\Phi(v,t)=\|v\|^{2}_{L^{2}}$. The infimum $\inf_{t\geq 0, \|h\|\geq r}\Phi(h,t)\rightarrow \infty$ as $r\rightarrow \infty$. We have
\begin{align}\label{Lyapunov}
\mathcal{L}_{t}\Phi(v,t)&=\frac{\partial}{\partial t}\Phi(v,t)+\frac{b^{2}}{2}\mbox{Tr}[\Phi''(v,t)v^{m}Qv^{m}]+(\Delta v, \Phi'(v,t))\nonumber\\
&\quad+(v-v^{3}, \Phi'(v,t))+\int_{Z}[\Phi(v+cv^{n}z,t)-\Phi(v,t)-(cv^{n}z, \Phi'(v,t))]\nu(dz)\nonumber\\
&\leq b^{2}\int_{D}q(x,x)v^{2m}(x)dx-2\int_{D}|\nabla v|^{2}dx+2\int_{D}(v^{2}-v^{4})dx\nonumber\\
&\quad+\int_{0}^{\infty}(cv^{n}z,cv^{n}z)\nu(dz)\nonumber\\
&\leq b^{2}q_{0}\|v\|^{2m}_{L^{2m}}+2\|v\|^{2}_{L^{2}}-2\|v\|^{4}_{L^{4}}+c^{2}\|v\|^{2n}_{L^{2n}}\int_{0}^{\infty}z^{2}\nu(dz).
\end{align}
By (\ref{alpha}) and (\ref{rpq}), we have
\begin{align}\label{rpqv}
\|v\|^{2m}_{L^{2m}}&\leq \|v\|^{2m\alpha}_{L^{2}}\|v\|^{2m(1-\alpha)}_{L^{4}}\nonumber\\
&\leq \varepsilon\|v\|^{\frac{2m(1-\alpha)}{1-m\alpha}}_{L^{4}}+\varepsilon^{-\frac{1-m\alpha}{m\alpha}}\|v\|^{2}_{L^{2}}\nonumber\\
&\leq \varepsilon\|v\|^{4}_{L^{4}}+\varepsilon^{-\frac{1-m\alpha}{m\alpha}}\|v\|^{2}_{L^{2}},
\end{align}
where $\alpha=\frac{2-m}{m}$. Similarly,
\begin{equation}\label{n}
\|v\|^{2n}_{L^{2n}}\leq \varepsilon \|v\|^{4}_{L^{4}}+\varepsilon^{-\frac{1-n\beta}{n\beta}}\|v\|^{2}_{L^{2}},
\end{equation}
where $\beta=\frac{2-n}{n}$. By (\ref{Lyapunov}), (\ref{rpqv}) and (\ref{n}),
\begin{align*}
\mathcal{L}_{t}\Phi(v,t)&=\big(b^{2}q_{0}\varepsilon+c^{2}\varepsilon\int_{0}^{\infty}z^{2}\nu(dz)-2\big)\|v\|^{4}_{L^{4}}\nonumber\\
&\quad+\big(b^{2}q_{0}\varepsilon^{-\frac{1-m\alpha}{m\alpha}}+c^{2}\varepsilon^{-\frac{1-n\beta}{n\beta}}\int_{0}^{\infty}z^{2}\nu(dz)+2\big)\|v\|^{2}_{L^{2}}.
\end{align*}
Choose $\varepsilon$  sufficiently small, we have
\begin{align*}
\mathcal{L}_{t}\Phi(v,t)\leq C\Phi(v,t).
\end{align*}
Therefore, by Theorem 5.1, the equation (\ref{AC}) has a global strong solution. \ \ $\Box$ \\


\end{document}